\documentclass[a4j]{article}
\usepackage{listings}
\usepackage[pdftex]{graphicx}
\usepackage{color}
\usepackage[]{multicol}
\usepackage{amsthm}
\usepackage{caption}
\usepackage{here}
\usepackage{mathrsfs}
\usepackage{subcaption}
\usepackage{amsmath, amssymb}
\usepackage{ascmac}
\usepackage{comment}
\usepackage{nccmath}

\newtheorem{theorem}{Theorem}

\newtheorem{corollary}{Corollary}
\newtheorem{lemma}{Lemma}

\def\vec#1{\mbox{\boldmath $#1$}}
\title{Expressing the largest eigenvalue of a singular beta $F$-matrix with heterogeneous hypergeometric functions}
\author{
	Koki Shimizu$^\text{a}$
	Hiroki Hashiguchi$^\text{a}$ 
	\vspace{1mm}
	\\
	\quad
	\begin{small}
	{
	\begin{tabular}{c}
	$^\text{a}$ Department of Applied Mathematics\\
	Tokyo University of Science\\
	1-3, Kagurazaka, Shinjuku-ku, Tokyo 162-8601, Japan
		\end{tabular}
	}
	\end{small}}
\date{}
\begin{document}
\maketitle
{\bf Abstract:}{
In this paper, the exact distribution of the largest eigenvalue of a singular random matrix for multivariate analysis of variance~(MANOVA) is discussed.
The key to developing the distribution theory of eigenvalues of a singular random matrix is to use heterogeneous hypergeometric functions with two matrix arguments. 
In this study, we define the singular beta $F$-matrix and extend the distributions of a nonsingular beta $F$-matrix to the singular case.
We also give the joint density of eigenvalues and the exact distribution of the largest eigenvalue in terms of heterogeneous hypergeometric functions. 
\section{Introduction}
The distribution of eigenvalues of an $F$-matrix plays an important role in multivariate analysis such as the test for equivalence of covariance matrices, MANOVA, and discriminant analysis. 
The distributions of an $F$-matrix are known to be equal to the distributions of a ratio of two Wishart matrices. 
Khatri~\cite{Khatri1972} derived the exact distributions of the largest and smallest eigenvalues of a nonsingular ratio of two real Wishart matrices using a finite series of Laguerre polynomials of matrix arguments.
Hashiguchi et al.~\cite{Hashiguchi2018} suggested an alternative derivation approach and conducted a numerical experiment using the holonomic gradient method~(HGM) for the hypergeometric functions of matrix arguments.
There are various approaches to deriving the approximate distributions of a nonsingular real $F$-matrix.
Johnstone's results \cite{Johnstone2008,Johnstone2009}, based on random matrix theory, showed that a Tracy-Widom distribution approximates the exact distribution for the largest eigenvalue. 
Matsubara and Hashiguchi~\cite{Matsubara2016} derived the Laplace approximation via F distributions for the nonsingular real F-matrix. 
Under the elliptical model, including the normal population, some results on the distributions of eigenvalues were discussed by Caro-Lopera et al.~\cite{Caro-Lopera2014} and Shinozaki et al.~\cite{Shinozaki2018}.
Caro-Lopera et al.~\cite{Caro-Lopera2014} presented the density of eigenvalues of a nonsingular ratio of two elliptical Wishart matrices for the moments of the modified likelihood ratio statistics.
D\'iaz-Garc\'ia and Guti\'errez-J\'aimez~\cite{Garcia2011} extended the nonsingular real Wishart distributions to the complex, quaternion and octonion cases under the normal population. 
These distributions are said to be nonsingular beta-Wishart distributions. 

In the case of a singular random matrix, Uhlig~\cite{Uhlig1994} derived the density of a singular real Wishart matrix and presented the Jacobian of the transformation to obtain the density of a singular real $F$-matrix as an open problem. 
The proof of Uhlig's result was given by D\'iaz-Garc\'ia and Guti\'errez-J\'aimez~\cite{Garcia1997}.
Chiani~\cite{Chiani2016} gave the distribution of the largest eigenvalue for a nonsingular or singular real F-matrix and provided an algorithm to compute the exact probabilities.
Shimizu and Hashiguchi~\cite{Shimizu2021} established the distribution theory of eigenvalues of a singular random matrix and derived the exact distributions of the largest eigenvalues of a singular beta-Wishart matrix.

In this paper, we discuss the distributions of eigenvalues of a singular beta $F$-matrix on real finite-dimensional algebra.
In Section~2, preliminary results and some notations are provided, which will be used throughout this paper. 
Furthermore, we introduce the heterogeneous hypergeometric functions of two matrix arguments, which were already defined by Shimizu and Hashiguchi~\cite{Shimizu2021}. 
These functions can be obtained using the integral formula for Jack polynomials over the Stiefel manifold.
In Section~3, the density functions of a singular beta $F$-matrix and the joint density functions of its eigenvalues are given.
We also show that the exact distribution of the largest eigenvalue can be expressed in terms of heterogeneous hypergeometric functions. 
This derivation is similar to that of Shimizu and Hashiguchi~\cite{Shimizu2021}.
Numerical computations on theoretical distributions are conducted using an algorithm presented by Hashiguchi et al.~\cite{Hashiguchi2000} for zonal polynomials.
In Section~4, we discuss the distribution of the largest eigenvalue of a ratio of two singular beta-Wishart matrices.

\section{Heterogeneous hypergeometric functions $_pF_q^{(\beta;m,n)}$}
In this section, we recall some notations that were provided by Shimizu and Hashiguchi~\cite{Shimizu2021}. 
We introduce the heterogeneous hypergeometric functions with two matrix arguments. 
These functions appear in the distributions of eigenvalues for a singular beta-Wishart matrix. 
Let $\mathbb{F}_\beta$ denote a real finite-dimensional division algebra such that
$\mathbb{F}_1 = \mathbb{R}$, $\mathbb{F}_2 =\mathbb{C}$, and $\mathbb{F}_4 = \mathbb{H}$ for $\beta=1, 2, 4$, where
$\mathbb{R}$ and $\mathbb{C}$ are the fields of real and complex numbers, respectively, and $\mathbb{H}$ is the quaternion division algebra over $\mathbb{R}$.
We restrict the parameter $\beta$ to values of $\beta = 1, 2$, and $4$, and denote $\mathbb{F}_\beta^{m \times n}$ 
as the set of all $m \times n$ matrices over $\mathbb{F}_\beta$, where $m \ge n$.
The conjugate transpose of $X \in \mathbb{F}_\beta^{m \times n}$ is written as $X^\ast=\overline{X^\top}$ and we say that $X$ is Hermitian if $X^\ast = X$.  
The set of all Hermitian matrices is denoted by 
$S^{\beta}(m)=\{ X \in \mathbb{F}^{m\times n}_\beta \mid X^\ast=X\}$. 
The eigenvalues of a Hermitian matrix are all real.  If the eigenvalues of $X \in S^\beta(m)$ are all positive, 
then we say that it is positive definite and write $X > 0$.
The exterior product $(dX)$ for $X \in \mathbb{F}_\beta^{m \times n}$ was defined by Mathai~\cite{Mathai1997} and D\'iaz-Garc\'ia and Guti\'errez-J\'aimez~\cite{Garcia2011}. 
We define the Stiefel manifold and the unitary group over $\mathbb{F}_\beta$ as
\begin{align*}
V^{\beta}_{n,m} = \{ H_1\in \mathbb{F}_\beta^{m\times n}\mid H_1^{\ast}H_1=I_n\},
\quad  U^{\beta}_m =\{ H \in \mathbb{F}_\beta^{m\times m} \mid H^{\ast}H = H H^{\ast} = I_m\} ,
\end{align*}
respectively.
If $\beta=1,2,4$, then $U^{\beta}_m$ are the real orthogonal group, unitary group, and symplectic group, respectively. 
The $\beta$-multivariate gamma function for $c \in \mathbb{F}_\beta$, $\Gamma_m^{\beta}(c)$,  is defined by
\begin{align*}
\Gamma_m^{\beta}(c) &= \int_{X > 0} |X|^{c - (m-1)\beta/2 } \mathrm{etr}(-X) (d X)
\\ 
 &= \pi^{\frac{m(m-1)\beta}{4}}\prod_{i=1}^{m}\Gamma\bigg\{c-\frac{(i-1)\beta}{2}\bigg\},
\end{align*}
where $\Re(c)>\frac{(m-1)\beta}{2}$, $|X|$ is the determinant of the matrix $X$ and $\mathrm{etr}(\cdot)=\exp(\mathrm{tr}(\cdot))$.
We define $(H_1^{\ast}dH_1)$ and $\mathrm{Vol}(V^{\beta}_{n,m})$ as
\begin{align*}
(H_1^{\ast}dH_1)=\bigwedge_{i=1}^{n}\bigwedge_{j=i+1}^{m}h_j^{\ast}dh_i,
\quad 
\mathrm{Vol}(V^{\beta}_{n,m})=\int_{H_1\in V^{\beta}_{n,m}}(H_1^{\ast}dH_1)=\frac{2^n\pi^{mn\beta/2}}{\Gamma^{\beta}_n(m\beta/2)} ,
\end{align*}
respectively, where 
$H_1\in V^{\beta}_{n,m}$ and $H = (H_1 \mid H_2)=(h_1, \dots, h_n \mid h_{n+1}, \dots, h_m) \in U_m^\beta$.
Another differential form, $(dH_1)$, is defined by
\begin{align*}
(dH_1)=\frac{(H_1^{\ast}dH_1)}{\mathrm{Vol}(V^{\beta}_{n,m})}=\frac{\Gamma^{\beta}_n(m\beta/2)}{2^n\pi^{mn\beta/2}}(H_1^{\ast}dH_1)
\end{align*}
and is normalized as 
$
\int_{H_1\in V^{\beta}_{n,m}}(dH_1)=1.
$

 For a positive integer $k$, let $\kappa=(\kappa_1,\dots,\kappa_m)$ denote a partition of $k$ with $\kappa_1\geq\cdots \geq \kappa_m\geq 0$ and $\kappa_1+\cdots +\kappa_m=k$. The set of all partitions with lengths not longer than $m$ is denoted by $P^k_{m}=\{\ \kappa=(\kappa_1,\dots,\kappa_m)\mid \kappa_1+\dots+\kappa_m=k, \kappa_1\geq \cdots \geq \kappa_m \geq 0 \}$. 
The $\beta$-generalized Pochhammer symbol of parameter $a>0$ is defined as 
\begin{align*}
(a)_\kappa^{\beta}=\prod_{i=1}^{m}\biggl(a-\frac{i-1}{2}\beta\biggl)_{\kappa_i},
\end{align*}
where $(a)_k=a(a+1)\cdots (a+k-1)$ and $(a)_0=1$.
The Jack polynomial $C_\kappa^{\beta}(X)$ is a symmetric polynomial in $x_1,\dots,x_m$; these are eigenvalues of $X\in S^\beta(m)$.
 See Stanley~\cite{Stanley1989} and Koev and Edelman~\cite{Edelman2006} for the relevant detailed properties.
If $\beta=1,2$, then the Jack polynomials are referred to as zonal polynomials and Shur polynomials, respectively. 
Li and Xue~\cite{Li2009} proposed zonal polynomials and hypergeometric functions of quaternion matrix arguments for $\beta=4$.
Shimizu and Hashiguchi~\cite{Shimizu2021} discussed the integral formula for Jack polynomials that differ from Theorem~3 in D\'iaz-Garc\'ia~\cite{Garcia2013}. 
This formula is needed to define the heterogeneous hypergeometric functions of two matrix arguments.

For $A\in S^{\beta}(m)$ and $B\in S^{\beta}(n)$, the integral formula for Jack polynomials over the Steifel manifold was given by Shimizu and Hashiguchi~\cite{Shimizu2021} as  
\begin{align}
\label{splitting:beta}
  \int_{H_1\in V^{\beta}_{n,m}}^{}C^{\beta}_{\kappa}(AH_1BH_1^\ast)(dH_1)&=\frac{C^{\beta}_{\kappa}(A)C^{\beta}_{\kappa}(B)}{C^{\beta}_{\kappa}(I_m)},
\end{align}
where $m\geq n$.

The hypergeometric functions of parameter $\beta>0$ are defined as 
\begin{align}
\label{def:hyper}
{_pF_q}^{(\beta;m)}({\mbox{\boldmath$\alpha$}};\mbox{\boldmath$\beta$};A)=\sum_{k=0}^{\infty}\sum_{\kappa \in P^k_m}\frac{(\alpha_1)^{\beta}_\kappa \cdots(\alpha_p)^{\beta}_\kappa}{(\beta_1)^{\beta}_\kappa \cdots(\beta_q)^{\beta}_\kappa}\frac{C^{\beta}_\kappa(A)}{k!},
\end{align} 
where $\mbox{\boldmath$\alpha$}=(\alpha_1,\ldots,\alpha_p)$, $\mbox{\boldmath$\beta$}=(\beta_1,\ldots,\beta_q)$.
The special case of (\ref{def:hyper}) can be represented as  ${_1F_0}^{(\beta;m)}(\alpha_1;A)=|I_m-A|^{-\alpha_1}$.
If $\beta=1$, we use ${_pF_q}$ instead of ${_pF_q}^{(1;m)}$.
Shimizu and Hashiguchi~\cite{Shimizu2021} defined the heterogeneous hypergeometric functions of two matrix arguments and provided some of their properties.  
Ratnarajah and Vaillancourt~\cite{Ratnarajah2005a, Ratnarajah2005b} used their functions to derive the joint density of the eigenvalues of a singular complex Wishart matrix.
From (\ref{splitting:beta}) and (\ref{def:hyper}), the heterogeneous hypergeometric functions is defined as follows.
\begin{align}
\nonumber
\label{def:hetero}
{_pF_q}^{(\beta;m,n)}(\mbox{\boldmath$\alpha$};\mbox{\boldmath$\beta$};A,B)&=\displaystyle \int_{H_1\in V^\beta_{n,m}}{_pF_q}^{(\beta;m)}(\mbox{\boldmath$\alpha$};\mbox{\boldmath$\beta$};AH_1BH_1^\ast)(dH_1)\\
&=\sum_{k=0}^{\infty}\sum_{\kappa \in P^k_n}\frac{(\alpha_1)^{\beta}_\kappa \cdots(\alpha_p)^{\beta}_\kappa}{(\beta_1)^{\beta}_\kappa \cdots(\beta_q)^{\beta}_\kappa}\frac{C^{\beta}_\kappa(A)C^{\beta}_\kappa(B)}{k!C^{\beta}_\kappa(I_m)}.
 \end{align}

\section{Exact distribution of eigenvalues of a singular beta $F$-matrix}
Suppose that an $m\times n$ beta-Gaussian random matrix $X$ is distributed as $X\sim$ $N^\beta_{m,n}(M,\Sigma \otimes I_n)$, 
where $M$ is the $m \times n$ mean matrix, $\Sigma > 0$, and $\otimes$ is the Kronecker product. This means that 
the column vectors of $X$ are an i.i.d. sample of size $n$ from $N^\beta_{m}({\vec{\mu}}, \Sigma)$, where $\vec{\mu}$ is the $m$-dimensional mean vector and $M=\vec{\mu}\bf{1}^\top$ and ${\bf{1}}=(1,\dots,1)^\top \in \mathbb{R}^n$. The density of $X$ is represented as 
\begin{align*}
 f(X)=\frac{1}{(2\pi \beta^{-1})^{mn\beta/2}|\Sigma|^{\ n\beta/2}}\mathrm{exp}\biggl(-\frac{\beta}{2}\mathrm{tr}(X-M)^{\ast}\Sigma^{-1}(X-M)\biggl).
\end{align*} 
Let $X\sim N^\beta_{m,n}(O,\Sigma \otimes I_n)$. 
If $n\geq m$, the random matrix $W=XX^{\ast}$ is called a nonsingular beta-Wishart matrix.
The distributions of a nonsingular real Wishart matrix have been studied by some authors. See Muirhead~\cite{Muirhead1982} and Gupta and Nagar~\cite{Gupta2000} for details.
Ratnarajah et al. \cite{Ratnarajah2005} derived the distribution of the largest and smallest eigenvalues for a nonsingular complex Wishart matrix. 
Their distributions were applied to the MIMO communication system. 
Kang and Alouini~\cite{Kang2003}  applied the exact distribution of the largest eigenvalue to the MIMO communication system and
Chiani et al.~\cite{Chiani2003} provided the exact expression for the characteristic function of MIMO system capacity.
Li and Xue~\cite{Li2009} discussed a quaternion case for a nonsingular Wishart matrix.
The density of a nonsingular beta-Wishart matrix that covers the real, complex, and quaternion Wishart matrices was given by D\'iaz-Garc\'ia and Guti\'errez-J\'aimez~\cite{Garcia2011} as
\begin{align*}
f(W)&=\frac{|\Sigma|^{-n\beta/2}}{(2\beta^{-1})^{ mn\beta/2}\Gamma^\beta_m(n\beta/2)}|W|^{(n-m+1)\beta/2-1}\mathrm{etr}\biggl(-\frac{\beta}{2}\Sigma^{-1}W\biggl).
\end{align*}
On the other hand, if $m>n$, the matrix $W$ of order $m$ was said to be a singular beta-Wishart matrix. 
A singular beta-Wishart matrix has only $n$ eigenvalues.
Shimizu and Hashiguchi~\cite{Shimizu2021} gave the density of a singular beta-Wishart matrix as follows.
\begin{align}
\label{betawishart}
 f(W)=\frac{\pi^{n(n-m)\beta/2}|\Sigma|^{- n\beta/2}}{(2\beta^{-1})^{mn\beta/2}\Gamma^{\beta}_n(n\beta/2)}(\mathrm{det}L_1)^{(n-m+1)\beta/2-1}\mathrm{etr}\biggl(-\frac{\beta}{2}\Sigma^{-1}W\biggl),
\end{align}
 where $W=G_1L_1G_1^{\ast}$,  $G_1\in V_{n,m}^{\beta}$, and $L_1$ is an $n \times n$ diagonal matrix.
 Throughout this paper, the notation $W^\beta_m(n,\Sigma)$ is referred to as either a nonsingular or singular case when $n\geq m$ or $n<m$, respectively.
\

For integers $p\geq m>n$, let $A\sim W^\beta_m(n,\Sigma_1)$ and $B\sim W^\beta_m(p,\Sigma_2)$, where $A$ and $B$ are independent.  
Put $B=T^\ast T$ where $T$ is an upper-triangular $m\times m$ matrix with the positive diagonal matrix. 
The singular beta $F$-matrix with rank $n$ is defined as $F=T^{-1}A(T^{-1})^\ast$.
The singular beta $F$-matrix has the same distribution as $AB^{-1}$.
The nonsingular parts of the spectral decomposition can be represented as $F=H_1QH_1^\ast$ where $H_1\in V^\beta_{n,m}$ and $Q=\mathrm{diag}(q_1,\dots,q_n)$ with $q_1>\cdots >q_n>0$.
The following lemma was first given by Uhlig~\cite{Uhlig1994} in a real case as a conjecture. D\'iaz-Garc\'ia and Guti\'errez-J\'aimez~\cite{Garcia1997} gave a proof of Uhlig's conjecture. D\'iaz-Garc\'ia and Guti\'errez-S\'anchez~\cite{Garcia2013andGutierrez-Sanchez} extended this result to complex, 
quaternion and octonion cases. 
\begin{lemma}
For $X,Y\in S^\beta(m)$ with rank $n<m$. Let $Y=T^{-1}X(T^{-1})^\ast$, where $T$ is a nonsingular $m\times m$ matrix. Let $X=G_1\Lambda_1G_1^\ast$ and $Y=H_1\Lambda_2H_1^\ast$, where $G_1$, $H_1\in V^\beta_{n,m}$ and $\Lambda_1, \Lambda_2$ are $n\times n$ diagonal matrices. Then we have
\begin{eqnarray*}
(dX)=|\Lambda_1|^{(m-n-1)\beta/2+1}|\Lambda_2|^{-(m-n-1)\beta/2-1}|T|^{n\beta}(dY).
\end{eqnarray*}
\end{lemma}

The following theorem represents the density of a singular beta $F$-matrix. 
 \begin{theorem}
Let $A\sim W^\beta_m(n,\Sigma_1)$ and $B\sim W^\beta_m(p,\Sigma_2)$, where $p\geq m>n$. Then the density of $F=B^{-1/2}AB^{-1/2}$ is given as 
\begin{align}
\frac{\pi^{n(n-m)\beta/2}|\Sigma|^{-n\beta/2}\Gamma^\beta_m\{(n+p)\beta/2\}}{\Gamma^\beta_n(n\beta/2)\Gamma^\beta_m(p\beta/2)}|Q|^{-(m-n-1)\beta/2-1}|I_m+\Sigma^{-1} F|^{-(n+p)\beta/2},
\end{align}
where $F=H_1QH_1^\ast$, $H_1\in V^\beta_{n,m}$, $Q=\mathrm{diag}(q_1,\dots,q_n)$ and $\Sigma=\Sigma_1\Sigma^{-1}_2$.
\begin{proof}
We can assume, without loss of generality, that $\Sigma_1=I_m$.
The density functions of $A$ and $B$ are given as 
\begin{align*}
 f(A)=\frac{\pi^{n(n-m)\beta/2}}{(2\beta^{-1})^{mn\beta/2}\Gamma^{\beta}_n(n\beta/2)}|L_1|^{(n-m+1)\beta/2-1}\mathrm{etr}\biggl(-\frac{\beta}{2}A\biggl),
\end{align*}
\text{and}\\
\begin{align*}
f(B)&=\frac{|\Sigma_2|^{-p\beta/2}}{(2\beta^{-1})^{mp\beta/2}\Gamma^\beta_m(p\beta/2)}|B|^{(p-m+1)\beta/2-1}\mathrm{etr}\biggl(-\frac{\beta}{2}\Sigma_2^{-1}B\biggl),
\end{align*}
respectively.
Then the joint density of $A$ and $B$ is given as 
\begin{align*}
f(A,B)&=\frac{\pi^{n(n-m)\beta/2}|\Sigma_2|^{-p\beta/2}}{(2\beta^{-1})^{m(n+p)\beta/2}\Gamma^\beta_n(n\beta/2)\Gamma^\beta_m(p\beta/2)}|L_1|^{(n-m-1)\beta/2-1}\\
&~~\times |B|^{(p-m+1)\beta/2-1}\mathrm{etr}\biggl\{-\frac{\beta}{2}(\Sigma_2^{-1}B+A)\biggl\}.
\end{align*}
From Lemma~1, the joint density of $F$ and $B$ is represented as 
\begin{align}
\label{ex01}
\nonumber
f(F,B)&=\frac{\pi^{n(n-m)\beta/2}|\Sigma_2|^{-p\beta/2}}{(2\beta^{-1})^{m(n+p)\beta/2}\Gamma^\beta_n(n\beta/2)\Gamma^\beta_m(p\beta/2)}|Q|^{-(m-n-1)\beta/2-1}\\
&~~\times |B|^{(n+p-m+1)\beta/2-1}\mathrm{etr}\bigg\{-\frac{\beta}{2}(\Sigma_2^{-1}+F)B\biggl\}.
\end{align}
Integrating (\ref{ex01}) respect to $B>0$ and $\Sigma_2^{-1}=\Sigma$, we get the density of $F$ as 
\begin{align*}
\frac{\pi^{n(n-m)\beta/2}|\Sigma|^{-n\beta/2}\Gamma^\beta_m\{(n+p)\beta/2\}}{\Gamma^\beta_n(n\beta/2)\Gamma^\beta_m(p\beta/2)}|Q|^{-(m-n-1)\beta/2-1}|I_m+\Sigma^{-1} F|^{-(n+p)\beta/2}.
\end{align*}

\end{proof}
\end{theorem}

\begin{theorem}
Under the same condition of $Theorem~3.1$, the joint density of eigenvalues $q_1,\dots,q_n$ of $F$ is given as
\begin{align}
\label{eigen:Fmatrix}
f(q_1,\dots,q_n)=C_1~|Q|^{(m-n+1)\beta/2-1}\prod_{i<j}^{n}(q_i-q_j)^\beta{_1F_0}^{(\beta;m,n)}\biggl(\frac{(p+n)\beta}{2};-\Sigma^{-1},Q\biggl),
\end{align}
where $Q=\mathrm{diag}(q_1,\dots,q_n)$, $C_1=\frac{\pi^{n^2\beta/2+r}|\Sigma|^{-n\beta/2}\Gamma^\beta_m\{(n+p)\beta/2\}}{\Gamma^\beta_n(n\beta/2)\Gamma^\beta_m(p\beta/2)\Gamma^\beta_n(m\beta/2)}$, and 
 \begin{align*}r = 
\begin{cases}
0 & \beta=1\\
-n\beta/2 & \beta=2,4.
\end{cases}
\end{align*}
\begin{proof}
The Jacobian of the transformation $F=H_1QH_1^\ast$ given by D\'iaz-Garc\'ia and Guti\'errez-S\'anchez~\cite{Garcia2013andGutierrez-Sanchez} is 
\begin{align}
\label{jacobian01}
(dF)=2^{-n}\pi^{r}\prod_{i=1}^{n}q_i^{(m-n)\beta}\prod_{i<j}^{n}(q_i-q_j)^{\beta}(dQ)\wedge(H_1^{\ast}dH_1). 
\end{align}
Using equation (\ref{jacobian01}) for the density of $F$, we have 
\begin{align*}
f(Q,H_1)&=\frac{2^{-n}\pi^{n(n-m)\beta /2+r}|\Sigma|^{-n\beta /2}\Gamma^\beta_m\{(n+p)\beta /2\}}{\Gamma^\beta_n(n\beta /2)\Gamma^\beta_m(p\beta /2)}\\
&~~\times|Q|^{(m-n+1)\beta/2-1}\prod_{i<j}^{n}(q_i-q_j)^{\beta}|I_m+\Sigma^{-1} H_1QH_1^\ast|^{-(n+p)\beta/2}.
\end{align*}
Moreover integrating $f(Q,H_1)$ with respect to $H_1$, the density of eigenvalues of $F$ is given as 
 \begin{align*}
f(q_1,\dots,q_n)&=C_1~|Q|^{(m-n+1)\beta /2-1}\prod_{i<j}^{n}(q_i-q_j)^\beta \\
&~~\times \int_{H_1\in V^\beta_{n,m}}|I_m+\Sigma^{-1} H_1QH_1^\ast|^{-(n+p)\beta /2}(dH_1).
\end{align*}
From (\ref{def:hetero}), the right hand side of the above equation can be evaluated as 
\begin{align*}
f(q_1,\dots,q_n)=C_1~|Q|^{(m-n+1)\beta/2-1}\prod_{i<j}^{n}(q_i-q_j)^\beta{_1F_0}^{(\beta;m,n)}\biggl(\frac{(p+n)\beta }{2};-\Sigma^{-1},Q\biggl).
\end{align*}
\end{proof}
\end{theorem}

The joint density $(\ref{eigen:Fmatrix})$ when $\Sigma_1=\Sigma_2$ is represented in Corollary~\ref{corollary:distuncorr}.
\begin{corollary}  \label{corollary:distuncorr}
Let $A\sim W^\beta_m(n,I_m)$ and $B\sim W^\beta_m(p,I_m)$, where $p\geq m>n$.
Then the joint density of eigenvalues $B^{-1/2}AB^{-1/2}$ is given as 
\begin{align*}
\frac{\pi^{n^2\beta/2+r}\Gamma^\beta_m\{(n+p)\beta/2\}}{\Gamma^\beta_n(n\beta/2)\Gamma^\beta_m(p\beta/2)\Gamma^\beta_n(m\beta/2)}|Q|^{(m-n+1)\beta/2-1}\prod_{i<j}^{n}(q_i-q_j)^\beta|I_n+Q|^{-(n+p)\beta/2}.
\end{align*}
\begin{proof}
Put $\Sigma=I_m$ for the density function (\ref{eigen:Fmatrix}), the functions ${_1F_0}^{(\beta;m,n)}$ can be represented as 
\begin{align*}
{_1F_0}^{(\beta;m,n)}\biggl(\frac{(p+n)\beta}{2};I_m,-Q\biggl)=|I_n+Q|^{-(n+p)\beta/2}.
\end{align*}
\end{proof}
\end{corollary}
The result for $\beta=1$ in Corollary~\ref{corollary:distuncorr} coincides with Theorem~4 (i) of D\'iaz-Garc\'ia and Guti\'errez-J\'aimez~\cite{Garcia1997} and the equation~(57), page~537 in Anderson~\cite{Anderson2003}. 
The next result was proposed by Sugiyama~\cite{Sugiyama1967} in the real case. Lemma~2 is required in order to integrate (\ref{eigen:Fmatrix}) with respect to $q_2, \dots, q_n$.
Shimizu and Hashiguchi~\cite{Shimizu2021} generalized this result for complex and quaternion cases.
   \begin{lemma}
  Let $X_1=\mathrm{diag}(1,x_2,\dots, x_n)$ and $X_2=\mathrm{diag}(x_2,\dots, x_n)$ with $x_2>\cdots >x_n>0$; then the following equation holds. 
  \begin{align*}
\int_{1>x_2>\cdots >x_n>0}|X_2|^{a-(n-1)\beta/2-1}C^{{\beta}}_\kappa(X_1)\prod_{i=2}^{n}(1-x_i)^\beta \prod_{i<j}(x_i-x_j)^\beta \prod_{i=2}^{n}dx_i\\
=(na+k)(\Gamma^{{\beta}}_n(n\beta/2)/\pi^{n^2\beta/2+r})\frac{\Gamma^{{\beta}}_n(a,\kappa)\Gamma^{{\beta}}_n\{(n-1)\beta/2+1\}C^\beta_\kappa(I_n)}{\Gamma^{{\beta}}_n\{a+(n-1)\beta/2+1,\kappa\}},
\end{align*}
where $\mathrm{\Re}(a)>(n-1)\beta/2$ and $\Gamma^{{\beta}}_n(\alpha,\kappa)=(\alpha)_\kappa \Gamma^\beta(\alpha)$.
\end{lemma}
\begin{theorem} \label{theo:ell1-dist}
Let $A\sim W^\beta_m(n,\Sigma_1)$ and $B\sim W^\beta_m(p,\Sigma_2)$, where $p\geq m>n$.
The distribution function of the largest eigenvalue $q_1$ of $B^{-1/2}AB^{-1/2}$ is given as 
\begin{align}
\label{maxeigen:beta}
C_2~|\Sigma|^{-n\beta/2}x^{{mn\beta/2}}
{_2F_1}^{(\beta;m,n)}\biggl(\frac{(n+p)\beta}{2};\frac{m\beta}{2};\frac{(m+n-1)\beta+1}{2};- x\Sigma^{-1},I_n\biggl),
\end{align}
where $C_2=\frac{\Gamma^\beta_m\left \{(n+p)\beta/2 \right \} \Gamma^{\beta}_n\left \{(n-1)\beta/2+1\right \}}{\Gamma^\beta_m\left \{p\beta/2 \right \}\Gamma^{\beta}_n\left\{(m+n-1)\beta/2+1\right\}}$.
\begin{proof}
We first start with the joint density (\ref{eigen:Fmatrix}).
\begin{align*}
&f(q_1,\dots,q_n)\\
&=C_1~|Q|^{(m-n+1)\beta/2-1}\prod_{i<j}^{n}(q_i-q_j)^\beta{_1F_0}^{(\beta;m,n)}\biggl(\frac{(p+n)\beta}{2};-\Sigma^{-1},Q\biggl)\\ 
&=C_1~|Q|^{(m-n+1)\beta/2-1}\prod_{i<j}^{n}(q_i-q_j)^\beta \sum_{k=0}^{\infty}\sum_{\kappa \in P^k_n}\frac{\{(p+n)\beta/2\}_\kappa~C^\beta_\kappa(-\Sigma^{-1})C^\beta_\kappa(Q)}{k!C^\beta_\kappa(I_m)}.
\end{align*}
Translating $q_i$ to $x_i=q_i/q_1$ for $i=2,\dots,n$ and using Lemma~2, the joint density of eigenvalues $f(q_1,\dots,q_n)$ can be evaluated as 
\begin{align*}
&f(q_1)=C_1~\int_{1>x_2>\cdots>x_n>0}|X_2|^{(m-n+1)\beta/2-1}C^\beta_\kappa(X_1) \prod_{i=2}^{n}(1-q_i)^\beta \prod_{2\leq i<j}^{n}(q_i-q_j)^\beta \\
&~~\times \sum_{k=0}^{\infty}\sum_{\kappa \in P^k_n}\frac{\{(p+n)\beta/2\}_\kappa~C^\beta_\kappa(-\Sigma^{-1})}{k!C^\beta_\kappa(I_m)} \\
&=C_1~(\Gamma^{{\beta}}_n(n\beta/2)/\pi^{n^2\beta/2+r})\frac{\Gamma^{{\beta}}_n(mn\beta/2,\kappa)\Gamma^{{\beta}}_n\left \{(n-1)\beta/2+1\right \}C^\beta_\kappa(I_n)}{\Gamma^{{\beta}}_n\left\{(m+n-1)\beta/2+1,\kappa\right\}}\\
&~~\times \sum_{k=0}^{\infty}\sum_{P^k_n}(mn\beta/2+k)q_1^{mn\beta/2+k-1}\frac{\{(p+n)\beta/2\}_\kappa~C^\beta_\kappa(-\Sigma^{-1})}{k!C^\beta_\kappa(I_m)} \\ 
&=\frac{\Gamma^\beta_m\left \{(n+p)\beta/2 \right \} \Gamma^{{\beta}}_n\left \{(n-1)\beta/2+1\right \}|\Sigma|^{-n\beta/2}}{\Gamma^\beta_m\left \{p\beta/2 \right \}\Gamma^{{\beta}}_n\left\{(m+n-1)\beta/2+1\right\}}\\  
&~~\times \sum_{k=0}^{\infty}\sum_{\kappa \in P^k_n}(mn\beta/2+k)q_1^{mn\beta/2+k-1}\frac{\{(p+n)\beta/2\}_\kappa\{m\beta/2\}_\kappa~C^\beta_\kappa(-\Sigma^{-1})C^\beta_\kappa(I_n)}{k!\left \{(m+n-1)\beta/2+1\right \}_\kappa C^\beta_\kappa(I_m)}. 
\end{align*}
Finally, integrating $f(q_1)$ with respect to $q_1$, we have the desired result (\ref{maxeigen:beta}).  
\end{proof}
\end{theorem}
\begin{corollary}
Under the same assumption of Theorem~\ref{theo:ell1-dist}, the distribution function of the largest eigenvalue $q_1$ of $B^{-1/2}AB^{-1/2}$ is represented as 
\begin{align}
\nonumber
\mathrm{Pr}&(q_1<x)\\ \nonumber
&=C_2~x^{{mn\beta/2}}\frac{|\Sigma|^{p\beta/2}}{|\Sigma+xI_m|^{(p+n)\beta/2}}\\  \label{positive:series}
&~~\times {_2F^{(\beta;m)}_1}\biggl(\frac{n\beta}{2}, \frac{(n+p)\beta}{2};\frac{(m+n-1)\beta}{2}+1;x(\Sigma+xI_m)^{-1}\biggl)\\ \nonumber
&=C_2~x^{{mn\beta/2}}|\Sigma+xI_m|^{-n\beta/2}\\ \label{finite:series}
&~~\times {_2F^{(\beta;m)}_1}\biggl(\frac{n\beta}{2}, \frac{(m-p-1)\beta}{2}+1;\frac{(m+n-1)\beta}{2}+1;x(\Sigma+xI_m)^{-1}\biggl) 
\end{align} 
\begin{proof}
This proof is similar way of Corollary~5 in Shimizu and Hashiguchi~\cite{Shimizu2021}.
If the length of a partition $\kappa$ is $m$, then we have $(\beta n / 2)^\beta_\kappa=0$, where $m > n$. Therefore the heterogeneous hypergeometric functions ${_2F_1}^{(\beta;m,n)}$ in (\ref{maxeigen:beta}) is represented as follows.
\begin{align}
\label{hyper:eigen}
 \nonumber
&{_2F_1}^{(\beta;m,n)}\biggl(\frac{(n+p)\beta}{2}, \frac{m\beta}{2};\frac{(m+n-1)\beta}{2}+1;- x\Sigma^{-1},I_n\biggl)\\  \nonumber
&=\sum_{k=0}^{\infty}\sum_{P^k_m}\frac{\{n\beta/2\}^\beta_\kappa \{(p+n)\beta/2\}^\beta_\kappa C^\beta_\kappa(-\Sigma^{-1})}{k!\{(m+n-1)\beta/2+1\}^\beta_\kappa}\\
&={_2F^{(\beta;m)}_1}\biggl(\frac{n\beta}{2}, \frac{(n+p)\beta}{2};\frac{(m+n-1)\beta}{2}+1;- x\Sigma^{-1}\biggl)
\end{align}
The $\beta$-Euler relations given in D\'iaz-Garc\'ia~\cite{Garcia2014} is 
\begin{align*}
_2F^{(\beta;m)}_1(a,b;c;X)
&=\mathrm{det}(I_m-X)^{-b}{_2F^{(\beta;m)}_1}(c-a,b;c;-X(I_m-X)^{-1})\\
&=\mathrm{det}(I_m-X)^{c-a-b}{_2F^{(\beta;m)}_1}(c-a,c-b;c;X).
\end{align*}
Using the above relationship for $_2F^{(\beta;m)}_1$ to (\ref{hyper:eigen}), we have the desired results. 
\end{proof}
\end{corollary}
We see that (\ref{positive:series}) is a series with positive terms.
On the other hand, we also see that (\ref{finite:series}) is a finite series if $r=(p-m+1)\beta/2-1$ is a positive integer.
\begin{corollary} 
Let $A\sim W^\beta_m(n,I_m)$ and $B\sim W^\beta_m(p,I_m)$, where $p\geq m>n$.
Then the distribution function of the largest eigenvalue $q_1$ of $B^{-1/2}AB^{-1/2}$ is given as
\begin{align}
\label{maxdist:indent}
\mathrm{Pr}&(q_1<x)=C_3~\biggl(\frac{x}{1+x}\biggl)^{mn/2}{_2F_1}\biggl(\frac{m-p+1}{2}, \frac{m}{2};\frac{n+m+1}{2};\frac{x}{1+x}I_n\biggl),
\end{align}
where $C_3=\frac{\Gamma_m\{(n+p)/2\}\Gamma_n\{(n+1)/2\}}{\Gamma_m\{p/2\}\Gamma_n\{(m+n+1)/2\}}$.
\begin{proof}
 If $\Sigma=I_m$ and $\beta=1$, the function $(\ref{maxeigen:beta})$ is given as 
\begin{align}
\label{maxeigen:identity}
\mathrm{Pr}(q_1<x)=C_3~
{_2F_1}\biggl(\frac{n+p}{2},\frac{m}{2};\frac{n+m+1}{2};-xI_n\biggl).
\end{align}
The Euler relations for ${_2F_1}$ is 
\begin{align}
\label{kummer}
{_2F_1}\biggl(a,b;c;X\biggl)=\mathrm{det}(I_n-X)^{-b}{_2F_1}\biggl(c-a,b;c;-X(I_n-X)^{-1}\biggl),
\end{align}
where $O<X<I_m$. From (\ref{kummer}), the function ${_2F_1}$ on the right hand side of (\ref{maxeigen:identity}) is represented as 
\begin{align*}
&{_2F_1}\biggl(\frac{n+p}{2},\frac{m}{2};\frac{n+m+1}{2};-xI_n\biggl)\\
&=(1+x)^{-mn/2}{_2F_1}\biggl(\frac{m-p+1}{2},\frac{m}{2};\frac{n+m+1}{2};\frac{x}{1+x}I_n\biggl).
\end{align*}
\end{proof}
\end{corollary}
Chiani~\cite{Chiani2016} provided the exact distribution of $q_1$ for $\beta = 1$ as a Pfaffian of a skew-symmetric matrix with a computable form.

\section{Numerical experiments}
In this section, we discuss the numerical experiments conducted for the distribution function (\ref{maxdist:indent}).
In general (\ref{maxdist:indent}) is an infinite series, but, if $r=(p-m-1)/2$ is a non-negative integer, it is a finite series represented by 
\begin{align}
\label{maxdist-ifinite:indent}
F(x)=C_3~\biggl(\frac{x}{1+x}\biggl)^{mn/2}\sum_{k=0}^{rn}\sum_{\kappa^\ast}{}\frac{\{(m-p+1)/2\}_\kappa\{m/2\}_\kappa~C_\kappa(\frac{x}{1+x}I_n)}{k!\left \{(m+n+1)/2\right \}_\kappa},
\end{align}
where $\sum_{\kappa^\ast}$ is the sum of all partitions of $k$ with $\kappa_1\leq r$.
The function (\ref{maxdist-ifinite:indent})  for $p=20, m=15$ and $n=3$ is given as 
\begin{align*}
\frac{128877}{8} \left(\frac{x}{x+1}\right)^{45/2} \left\{\frac{260 x^6}{1083 (x+1)^6}-\frac{104 x^5}{57
   (x+1)^5}+\frac{110 x^4}{19 (x+1)^4}\right.\\
 \left.-\frac{500 x^3}{51 (x+1)^3}+\frac{3725 x^2}{399 (x+1)^2}-\frac{90 x}{19
   (x+1)}+1
   \right\}.
   \end{align*}
We find that $r n=6$ in the series (\ref{maxdist-ifinite:indent}) and that the finite series (\ref{maxdist-ifinite:indent}) can be calculated using small terms when the parameters $p$ and $m$ are the value of the close. 
Chiani~\cite{Chiani2016} provided an algorithm to compute the exact distribution even if $r$ is a non-negative integer.
We denote the exact distribution by $F_0$ based on the algorithm.
In Table~1, the percentile points of the functions $F(x)$ and $F_0(x)$ are shown.

 \begin{table}[H]
\begin{center}
\caption{Percentile points of $q_1$ of $F^1_m(I_m,2,10)$} \label{table:pp-l1}
\begin{tabular}{c}
            \begin{minipage}{0.3\hsize}
        \begin{center}
         \captionsetup{labelformat=empty,labelsep=none}
          \subcaption{$p=20, m=5, n=4$}
          \label{table1}
               
           \begin{tabular}{lcc}
\hline
$\alpha$&${{F_0}^{\small{-1}}}(\alpha)$ &$F^{-1}(\alpha)$ \\
\noalign{\smallskip}\hline\noalign{\smallskip}
0.01	 &~0.27& 	0.27\\
0.05	 &~0.37& 	0.37\\
0.50	 &~0.82& 	0.82\\
0.95	 &~1.86& 	1.86\\
0.99  &~2.67& 	2.67\\ 
\noalign{\smallskip}\hline
\end{tabular}
        \end{center}
      \end{minipage}
      \begin{minipage}{0.3\hsize}
        \begin{center}
        \captionsetup{labelformat=empty,labelsep=none}
          \subcaption{$p=20, m=15, n=4$}
         \begin{tabular}{lcc}
\hline
$\alpha$&${{F_0}^{\small{-1}}}(\alpha)$ &$F^{-1}(\alpha)$  \\
\noalign{\smallskip}\hline\noalign{\smallskip}
0.01	 &~2.33& 	2.33\\
0.05	 &~3.16& 	3.16\\
0.50	 &~7.51& 	7.51\\
0.95	 &~24.4& 	24.4\\
0.99  &~45.6& 	45.6\\
\noalign{\smallskip}\hline
\end{tabular}
        \end{center}
      \end{minipage}
    \end{tabular}
  \end{center}
\end{table}
We test the hypothesis of the equality of mean vectors of some multivariate normal populations in MANOVA. 
We use Wilks Likelihood Ratio, Lawley Hotelling Trance, Bartlett-Nanda-Pillai, and Roy's test statistics relative to eigenvalues of some nonsingular or singular beta $F$-matrix.
In a two-sample problem, these four statistics and the well-known Hotelling $T^2$ statistic yield equivalent results because the singular beta $F$-matrix has one nonzero eigenvalue. 
The distribution (\ref{maxdist:indent}) of the Roy's test statistic can be available if the number of groups is less than or equal to the number of variables.
Rencher and Christensen~\cite{Rencher2012} provided the data in Table~6.12 on eight trees from each rootstock. 
Each tree considers four variables. 
We test the equality of the four mean vectors from the first four groups of rootstocks data using the distribution (\ref{maxdist:indent}) of the Roy's test statistic.
The nonzero eigenvalues of the singular real $F$-matrix for rootstocks data were $2.73, 0.54$, and $0.033$. 
The largest eigenvalue $2.73$ accounts for $2.73 /(2.73+0.54+0.033)=0.827$ of the sum of the nonzero eigenvalues.
We compute the approximately truncated distribution using zonal polynomials with degrees at most $30$ for the distribution (\ref{maxdist:indent}) with $n=(4-1)$, $m=4$ and $p=4(8-1)$. 
The 95th percentile point is $0.73$, which is consistent with that from using the exact calculation algorithm of Chiani~\cite{Chiani2016}.
Because the $95$ percentile point of the function is $0.763$ that is less than $2.73$, we reject the hypothesis that the four mean vectors are equal.

We discuss the distribution of the ratio of two Wishart matrices when they are both singular. 
In this case, their distributions have often used the test of equalities for covariance matrices.
Let $A\sim W^1_p(m,I_p)$ and $B\sim W^1_p(n,I_p)$ where $A$ and $B$ are independent.  
Assuming $p>m>n$, the product $AB^{-1}$ is not defined. 
Srivastava~\cite{Srivastava2007} defined the product $AB^+$ instead of $AB^{-1}$, where $A^{+}$ is a Moore-Penrose inverse matrix of $A$. 
The density of eigenvalues of $AB^{+}$ is approximately equivalent to $\frac{1}{p}U$ for a large $m$ where $U$ is distributed as a singular real Wishart distribution $W^1_m(n,I_m)$. 
Furthermore, the distributions of the eigenvalues of $AB^{+}$ and $UW^{-1}$ are also equivalent, where $W\sim W^1_m(p,I_m)$. 
See Srivastava~\cite{Srivastava2007} and Grinek~\cite{Grinek2019} for details. 
We extend this result to the complex and quaternion cases.
The next theorem implies that the exact distribution of the largest eigenvalue of $AB^+$ is equivalent to the distribution of (\ref{maxdist:indent}).
\begin{theorem}
Let $A\sim W^\beta_p(n,I_p)$ and $B\sim W^\beta_p(m,I_p)$ where $A$ and $B$ are independent. If $p>m>n$, then the distribution of the largest eigenvalue $q_1(A B^+)$ is equivalent to the largest eigenvalue $q_1(UW^{-1})$,
where $U\sim W^\beta_m(n,I_m)$, $W\sim W^\beta_m(p,I_m)$.
\begin{proof}
The singular beta Wishart matrix $A$ and $B$ can be written as $A=ZZ^\ast$ and $B=H_1LH_1^\ast$ where $Z\sim N^\beta_{p,n}(O,I_p\otimes I_n)$, $L$ is the $m\times m$ diagonal matrix, and $H_1\in V^\beta_{m,p}$. Then 
\begin{align*}
q_1(AB^+)&=q_1(Z^\ast B^+Z)\\
&=q_1(Z^\ast H_1L^{-1}H_1^\ast Z),
\end{align*}
where $V=H_1^\ast Z$ is distributed as matrix variate beta normal distributions $N^\beta_{m,n}(O,I_m\otimes I_n)$.
We consider the nonsingular random matrix $H L^{-1}H^\ast$, where $H\in U^\beta_m$.
Then we have
\begin{align*}
q_1(AB^+)&=q_1(V^\ast L^{-1}V)\\
&=q_1(V^\ast H^\ast H L^{-1}H^\ast HV)\\
&=q_1(HV V^\ast H^\ast H L^{-1}H^\ast)\\
&=q_1(UW^{-1}).
\end{align*}
\end{proof}
\end{theorem}
The exact computation of the largest eigenvalue of $AB^+$ using distribution (\ref{maxdist-ifinite:indent}) for $p=10, m=5$ and $n=3$ is represented as 
\begin{align}\nonumber
\label{dist:(10,5,3)}
\frac{693}{4} \left(\frac{x}{x+1}\right)^{15/2} \left\{\frac{5 x^6}{198 (x+1)^6}-\frac{3 x^5}{11 (x+1)^5}+
\frac{85x^4}{66 (x+1)^4}\right.\\
\left.-\frac{250 x^3}{77 (x+1)^3}+\frac{50 x^2}{11 (x+1)^2}-\frac{10 x}{3 (x+1)}+1\right\}.
\end{align}
 The 95 percentage point of (\ref{dist:(10,5,3)}) is 7.63.
We obtain the upper $5\%$ probability of the exact distribution given in Chiani~\cite{Chiani2016} for the largest eigenvalue of $AB^+$, which is 0.050. 
Finally, we calculate the exact distribution (\ref{finite:series}) that is also represented by a finite series for any $\Sigma>0$.
Fig~\ref{fig:graph01} shows the graph of (\ref{finite:series}) with parameters $p=10$, $n=2$ and $\Sigma=(1/3,1/2,1)$.
\begin{figure}[H]
 \begin{center}
     \includegraphics[width=7cm]{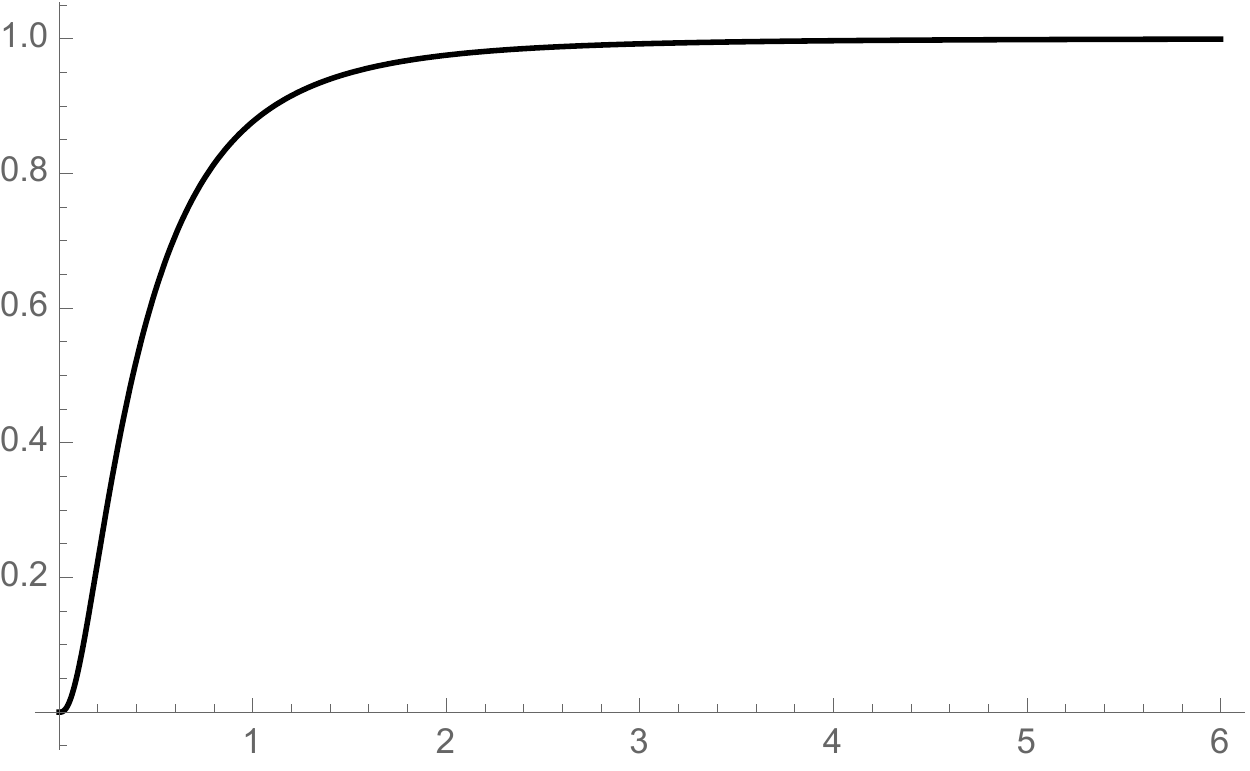}
     \rlap{\raisebox{1.0ex}{\kern-0.0em{\small $x$}}}%
\rlap{\raisebox{30.0ex}{\kern-21.0em{\small $\Pr(q_1 < x) $}}}%
        \caption{$p=10$, $n=2$, $\Sigma=(1/3,1/2,1)$} \label{fig:graph01}
        \end{center}
  \label{fig:one}
\end{figure}

\section{Conclusion}
In this study, we discussed the exact distribution of the largest eigenvalue of a singular random matrix. 
The exact distribution of the largest eigenvalue of a singular beta $F$-matrix was derived in terms of heterogeneous hypergeometric functions.
Numerical experiments were performed for the theoretical distributions of (\ref{finite:series}) and (\ref{maxdist:indent}).
We also considered the distribution of the largest eigenvalue of the ratio of two singular beta-Wishart matrices. This distribution for $\beta=1$ could be reduced in the form of (\ref{maxdist:indent}).

\end{document}